\documentclass[a4paper, 11pt]{amsart}
\usepackage{amscd,amsthm,amsfonts,latexsym,amssymb}
\usepackage[dvips]{graphicx}

\usepackage{calc,curves,ebezier,epic,eepic,graphicx,multiply,rotating}
\usepackage{tikz-cd} 
\usepackage{tensor}
\newtheorem{theorem}{Theorem} [section]

\newtheorem{example}[theorem]{Example}

\newtheorem{remark}[theorem]{Remark}

\renewcommand{\phi}{\varphi}

\newcommand{\trace}{\operatorname{trace}}

\begin{document}

\title{Almost complex surfaces in the nearly K\"ahler flag manifold}
\author{Kamil Cwiklinski and Luc Vrancken}

\address{Kamil Cwiklinski, Service de Physique de l'Univers, Champs et Gravitation, Universit\'{e} de Mons, 20 Place du Parc, 7000 Mons, Belgium
}
\email{kamil.cwiklinski@umons.ac.be}

\address{Luc Vrancken, LMI-Laboratoire de Math\'ematiques pour l'Ingénieur \\
Universit\'e Polytechnique Hauts-de-France \\
Campus du Mont Houy  \\
59313 Valenciennes Cedex 9, France}

\email{luc.vrancken@uphf.fr}

\address{Luc Vrancken, Katholieke Universiteit Leuven, Departement Wiskunde \\ Celestijnenlaan 200 B 3001 Leuven Belgium}

\begin{abstract} 
We study and classify almost complex totally geodesic submanifolds of the nearly K\"ahler flag manifold $F_{1,2}(\mathbb C^3)$, and of its semi-Riemannian counterpart. We also develop a structural approach to the nearly K\"ahler flag manifold $F_{1,2}(\mathbb C^3)$, expressing for example the curvature tensor in terms of the nearly K\"ahler structure $J$ and the three canonical orthogonal complex structures. 
\end{abstract}



\maketitle

\thispagestyle{empty}


\section{Introduction}   

Let $(\tilde M,<.,.>)$ be a (pseudo)-Riemannian manifold equipped with an almost complex structure $J$ which is compatible with the (pseudo)-Riemannian metric. We call $\tilde M$ a nearly K\"ahler manifold if and only if for any tangent vector field $X$ on $\tilde M$ we have  
$$(\tilde \nabla_X J)X = \tilde \nabla_X (JX)-J (\tilde \nabla_X X)= 0,$$
where $\tilde \nabla$ is the Levi-Civita connection associated to the nearly K\"ahler metric. 
This is equivalent to stating that the tensor 
$$G(X,Y)= (\tilde \nabla_XJ)Y$$
is a skew-symmetric tensor. The study of nearly K\"ahler manifolds and their submanifolds started with the work of Gray around 1970, see for example the references 
\cite{gray1}, \cite{gray2}, \cite{gray3}, \cite{gray4}.

Note that it is well known that a $2$ (or $4$) dimensional nearly K\"ahler manifold is always K\"ahler. Therefore the study of nearly K\"ahler manifolds is only worthwhile starting from dimension $6$. Perhaps surprisingly not that many $6$-dimensional nearly K\"ahler manifolds are known. Traditionally one has the following examples:
\begin{enumerate}
\item the $6$-dimensional sphere with its standard metric and an almost complex structures defined using the Cayley numbers;
\item the manifold $S^3 \times S^3$ not equipped with the product metric but with a metric introduced by a canonical submersion from $S^3 \times S^3 \times S^3$ to $S^3 \times S^3$;
\item the space of flags $F_{1,2}(\mathbb C^3)$;
\item the complex projective space $\mathbb CP^3$
\end{enumerate}
as well as their pseudo-Riemannian counterparts. It was shown by Butruille that in the Riemannian case, these are all the homogeneous $6$-dimensional nearly K\"ahler spaces. A similar classification result is not known in the pseudo-Riemannian case. Recently Foscolo and Haskins (see \cite{foshas}) discovered new examples of non-homogeneous complete structures on $S^6$ and $S^3 \times S^3$.  

In this paper we are in particular interested in the nearly K\"ahler structure on the space of flags $F_{1,2}(\mathbb C^3)$, which as a homogeneous space can be represented as $ SU(3) / (U(1) \times U(1)) $ and its pseudo-Riemannian counterpart $SU(2,1)/(U(1) \times U(1))$. A structural approach to these spaces is given in Section 2. The approach we present here is different from the one in \cite{loubeau}.

In the final section we will then investigate totally geodesic almost complex submanifolds. A submanifold $M$ is called almost complex if $J$ preserves the tangent space. It is well known that such a submanifold is even-dimensional and that a $6$-dimensional (pseudo-) nearly K\"ahler manifold does not admit any $4$-dimensional almost complex submanifolds. Therefore we can restrict ourselves to the $2$-dimensional case and we obtain a complete classification of the totally geodesic surfaces of $F_{1,2}(\mathbb C^3)$ and its pseudo-Riemannian counterpart. 

More precisely we prove the following theorems:
\begin{theorem} Let $M$ be a totally geodesic almost complex surface in the flag manifold $F_{1,2}(\mathbb C^3)$. Then either:
\begin{enumerate}
\item $K=4$ and $M$ is locally congruent to an open part of the immersion described in Example 3.1 which is a sphere with radius $\tfrac 12$;
\item $K=1$  and $M$ is locally congruent to an open part of the immersion described in Example 3.2 which is a sphere with radius $1$;
\item $K=0$ and $M$ is locally congruent to an open part of the immersion described in Example 3.3 which is a flat torus.
\end{enumerate}
\end{theorem}

Similarly in the pseudo-Riemannian case we obtain:
\begin{theorem} Let $M$ be a totally geodesic almost complex surface in the pseudo-Riemannian nearly K\"ahler manifold $\frac{SU(2,1)}{U(1) \times U(1)}$
Then either:
\begin{enumerate}
\item $K=4$ and $M$ is locally congruent to an open part of the immersion described in Example 3.4 which is a sphere with radius $\tfrac 12$;
\item $K=4$  and $M$ is locally congruent to an open part of the immersion described in Example 3.5 which is anti-isometric to a hyperbolic plane with curvature $-4$;
\item $K=1$ and $M$ is locally congruent to an open part of the immersion described in Example 3.6 which is anti-isometric to the hyperbolic plane with curvature $-1$. 
\end{enumerate}
\end{theorem}

\section{The nearly K\"ahler structure on $F_{1,2}(\mathbb C^3)$}

We now describe the nearly K\"ahler structure on  $F_{1,2}(\mathbb C^3)$, which as a homogeneous space can be represented by $\frac{SU(3)}{U(1) \times U(1)}$. We will use as a starting point the description of the structure described for example in \cite{butruille}, \cite{reinierstorm}, \cite{apostolovgrantcharovivanov}, but we will adapt it using a tensorial approach in order to obtain an explicit tensorial expression for the curvature tensor. 

Recall that $SU(3) =\{g \in GL(3,\mathbb C) \vert \bar g^t = g^{-1} \quad \text{and} \det g =1\}$ and that its Lie algebra $\mathfrak{su(3)} = \{ A \in \mathbb C^{3\times 3} \vert \bar A^t =-A \quad \text{and} \quad \trace A=0\}$. We consider as a basis of $\mathfrak{su(3)}$ the following vector fields:
\begin{alignat*}{2}
&h_1=\begin{pmatrix}
-i&0&0\\0&0&0\\0&0&i
\end{pmatrix}\qquad
&h_2=\begin{pmatrix}
\tfrac{i}{\sqrt{3}}&0&0\\0&\tfrac{-2i}{\sqrt{3}}&0\\0&0&\tfrac{i}{\sqrt{3}}
\end{pmatrix}\\
&m_1=\begin{pmatrix}
0&-1&0\\1&0&0\\0&0&0
\end{pmatrix}\qquad
&m_2=\begin{pmatrix}
0&0&0\\0&0&-1\\0&1&0
\end{pmatrix}\\
&m_3=\begin{pmatrix}
0&0&-1\\0&0&0\\1&0&0
\end{pmatrix}\qquad
&m_4=\begin{pmatrix}
0&i&0\\i&0&0\\0&0&0
\end{pmatrix}\\
&m_5=\begin{pmatrix}
0&0&0\\0&0&i\\0&i&0
\end{pmatrix}\qquad
&m_6=\begin{pmatrix}
0&0&i\\0&0&0\\i&0&0
\end{pmatrix}
\end{alignat*}
We denote by  $\mathfrak{h}$, resp.  $\mathfrak{m_1}$, $\mathfrak{m_2}$, $\mathfrak{m_3}$, the vector spaces spanned by $\{h_1,h_2\}$, resp. $\{m_1,m_4\}$, $\{m_2,m_5\}$, $\{m_3,m_6\}$. 

It is immediately clear that  $\mathfrak{h}$, $\mathfrak{h} \oplus  \mathfrak{m_1}$, $\mathfrak{h} \oplus  \mathfrak{m_2}$ and $\mathfrak{h} \oplus  \mathfrak{m_3}$ are Lie subalgebras of $\mathfrak{su(3)}$. We denote the corresponding Lie subgroups of $SU(3)$ by $H= U(1) \times U(1)$, $G_1=SU(2) \times U(1)$, $G_2=SU(2) \times U(1)$ and $G_3=SU(2) \times U(1)$. These are all closed subgroups of $SU(3)$ and as such there is unique way to make the projection on the quotient manifold a submersion. We call the quotient manifold $\frac{SU(3)}{U(1)\times U(1)}$ the flag manifold $F_{1,2}(\mathbb C^3)$.  We will denote the projection by $\pi$ and we identify the tangent space at $\pi(I)$ with 
$\mathfrak{m}_1 \oplus \mathfrak{m}_2 \oplus \mathfrak{m}_3$.  

As $SU(3)$ is a Lie subgroup of $GL(3,\mathbb C)$ its Killing form $B$ is given by $B(X,Y)= Re (\trace \bar X^t Y)$, so that the Killing form is proportional to the trace form. This allows us to define the metric $g$ on $SU(3)$ by $g_p(X,Y)=\frac{1}{2} \text{trace} ( \overline{X}^t Y)$. Using this metric, it is easy to verify that the previously constructed basis is an orthonormal basis. Since $g$ is a bi-invariant metric, the associated connection $D$ and its curvature tensor $R^D$ are respectively given by
\begin{align*}
&D_X Y = \tfrac 12 [X,Y]\\
&R^D(X,Y)Z = \tfrac 14 [Z, [X,Y]]
\end{align*}

Taking as complementary space $\mathfrak{m}= \mathfrak{m_1} \oplus \mathfrak{m_2} \oplus \mathfrak{m_3}$, we see in a straightforward way that $\mathfrak{m}$ is $Ad(H)$-invariant, i.e. for every $h \in H$ and every $X \in \mathfrak{m}$  we have that $Ad_h(X)=hXh^{-1} \in \mathfrak{m}$. This means that the flag manifold is a reductive homogeneous space. 
If we now parametrize $H$ by 
$$\left(
\begin{array}{ccc}
 e^{\frac{1}{3} i \left(\sqrt{3} s-3 t\right)} & 0 & 0 \\
 0 & e^{-\frac{2 i s}{\sqrt{3}}} & 0 \\
 0 & 0 & e^{\frac{1}{3} i \left(\sqrt{3} s+3 t\right)} \\
\end{array}
\right)
$$
we see that 
\begin{align*}
&Ad_h(m_1)= \left(
\begin{array}{ccc}
 0 & -e^{i \left(\sqrt{3} s-t\right)} & 0 \\
 e^{-i \left(\sqrt{3} s-t\right)} & 0 & 0 \\
 0 & 0 & 0 \\
\end{array}
\right)=\cos \left(\sqrt{3} s-t\right) m_1-\sin \left(\sqrt{3} s-t\right) m_4 \in \mathfrak{m_1}\\
&Ad_h(m_4)=\left(
\begin{array}{ccc}
 0 & i e^{i \left(\sqrt{3} s-t\right)} & 0 \\
 i e^{-i \left(\sqrt{3} s-t\right)} & 0 & 0 \\
 0 & 0 & 0 \\
\end{array}
\right)=\sin \left(\sqrt{3} s-t\right) m_1+\cos \left(\sqrt{3} s-t\right) m_4 \in \mathfrak{m_1}
\end{align*}
Similar results holds for $\mathfrak{m_2}$ and $\mathfrak{m_3}$. As the flag manifold is a reductive homogeneous space, there is a correspondence between $Ad(H)$-invariant tensors on $\mathfrak{m}$ and $SU(3)$-invariant tensors on $F_{1,2}$. So we see that
\begin{enumerate}
\item the flag manifold admits three rank $2$ distributions which are preserved by $SU(3)$. We will denote these distributions by $V_1$, $V_2$, $V_3$. 
\item on each of these distributions we have a natural metric (determined up to a scalar factor) and a natural almost complex structure (determined up to sign). The factors of the metric induced by the bi-invariant metric on $SU(3)$ are $1$. Taking these factors turns the projection of $SU(3) \rightarrow F_{1,2}$ in a Riemannian submersion. 
\end{enumerate}

So we see that we actually can define a $3$-parameter family of metrics on $F_{1,2}$ which are preserved by $SU(3)$ by
$$g_{\lambda_1,\lambda_2,\lambda_3}(X,Y)=\lambda_1 <X_1,Y_1>_1 +\lambda_2 <X_2,Y_2>_2 +\lambda_3 <X_3,Y_3>_3,$$
where $\lambda_1,\lambda_2, \lambda_3$ are positive numbers, $<.,.>_1$, $<.,.>_2$ and $<.,.>_3$ are the canonical metrics on $V_1$, $V_2$, $V_3$. It is clear that the metric corresponding to the Riemannian submersion, which is the main metric which we will consider in this section, corresponds to the choice of $\lambda_1=\lambda_2=\lambda_3=1$. 

We can also immdiately define $4$ almost complex structures determined respectively by
\begin{alignat}{3}
&J(m_1)= m_4 \qquad &J(m_2)=m_5 \qquad &J(m_3)=-m_6\\
&J_1(m_1)= m_4 \qquad &J_1(m_2)=-m_5 \qquad &J_1(m_3)=m_6\\
&J_2(m_1)= -m_4 \qquad &J_2(m_2)=-m_5 \qquad &J_2(m_3)=-m_6\\
&J_3(m_1)= -m_4 \qquad &J_3(m_2)=m_5 \qquad &J_3(m_3)=m_6.
\end{alignat}

It is clear from the definition that these almost complex structures commute and that they are compatible with respect to any metric of our $3$-parameter family of metrics. 
Moreover we have that
\begin{enumerate}
\item $J=J_1$ on $V_1$ and $J=-J_1$ on $V_2 \oplus V_3$,
\item $J=J_2$ on $V_3$ and $J=-J_2$ on $V_1 \oplus V_2$,
\item $J=J_3$ on $V_2$ and $J=-J_3$ on $V_1 \oplus V_3$.
\end{enumerate}
It is also immediate that $J=-(J_1+J_2+J_3)$ and $J=-J_1J_2J_3$.

Note that one can check (see for example \cite{apostolovgrantcharovivanov}) that $J_1$, $J_2$ and $J_3$ are integrable almost complex structures
on $F_{1,2}$ and that moreover $(g_{2,1,1},J_1)$, $(g_{1,2,1},J_2)$ and $(g_{1,1,2},J_3)$ are K\"ahler-Einstein structures of non-negative sectional curvature on $F_{1,2}$. 

However here in this paper we will consider the structure on $F_{1,2}$ determined by the metric $g=g_{1,1,1}$ which is the metric induced by the submersion together with the almost complex structure $J$. 

We denote the Levi-Civita connection on $F_{1,2}$ by $\bar \nabla$. We identify vectors in the Lie algebra in the usual way with Killing vector fields. Thus, writing $ o = \pi(I)$, we let $d\pi (X) = X_o^{*}$ where $X^{*}$ is the Killing vector field. Thus we can write $(\bar \nabla_{X^{*}} Y^{*})_o = - \frac{1}{2} [X, Y]_{\mathfrak{m}}$ with $X, Y \in \mathfrak{g} = \mathfrak{su} (3)$, see also \cite{arvanitoyeorgos, besse}.
This way we get
\begin{alignat}{3}
&\bar\nabla_{m_1} m_2 = \tfrac 12 m_3&\quad &\bar\nabla_{m_2} m_3 = \tfrac 12 m_1\quad &\bar \nabla_{m_3} m_1 = \tfrac 12 m_2\\
&\bar\nabla_{m_1} m_3 = - \tfrac 12 m_2&\quad &\bar\nabla_{m_2} m_1 = - \tfrac 12 m_3\quad &\bar \nabla_{m_3} m_2 = - \tfrac 12 m_1\\
&\bar\nabla_{m_1} m_5 =\tfrac 12 m_6&\quad &\bar\nabla_{m_2} m_6 =\tfrac 12 m_4\quad &\bar \nabla_{m_3} m_4 = - \tfrac 12 m_5\\
&\bar\nabla_{m_1} m_6 = - \tfrac 12 m_5&\quad &\bar\nabla_{m_2} m_4 = - \tfrac 12 m_6\quad &\bar \nabla_{m_3} m_5 = \tfrac 12 m_4\\
&\bar\nabla_{m_4} m_2 = \tfrac 12 m_6&\quad &\bar\nabla_{m_5} m_3 = - \tfrac 12 m_4\quad &\bar \nabla_{m_6} m_1 = \tfrac 12 m_5\\
&\bar\nabla_{m_4} m_3 = \tfrac 12 m_5&\quad &\bar\nabla_{m_5} m_1 = - \tfrac 12 m_6\quad &\bar \nabla_{m_6} m_2 = -\tfrac 12 m_4\\
&\bar\nabla_{m_4} m_5 =- \tfrac 12 m_3&\quad &\bar\nabla_{m_5} m_6 =  \tfrac 12 m_1\quad &\bar \nabla_{m_6} m_4 =  \tfrac 12 m_2\\
&\bar\nabla_{m_4} m_6 = -\tfrac 12 m_2&\quad &\bar\nabla_{m_5} m_4 =  \tfrac 12 m_3\quad &\bar \nabla_{m_6} m_5 = - \tfrac 12 m_1,
\end{alignat}
and all other coefficients at that point vanish. Given that the flag manifold is a homogeneous space, and the structure $J$ is preserved by isometries, a straightforward computation now shows that $(F_{1,2}, g,J)$ is a nearly K\"ahler space. This is shown by evaluating the tensor $G(X,Y) =( \bar \nabla_X  J ) (Y)$ on the basis vector fields $m_i$ and calculating that  $G(m_1,m_2) = m_6 = - G(m_2, m_1)$ etc. From this, using the properties of $G$ for a nearly K\"ahler space all the other coefficients follow in a straightforward way. We also deduce that the constant $\alpha$  associated to a $6$-dimensional nearly K\"ahler manifold for our choice of metric equals $1$. 

As $F_{1,2}$ is a reductive homogeneous space, we see from the proof of Proposition  3.4 of \cite{kobayashinomizu} that the curvature tensor $\bar R$ of $\bar \nabla$ is given by
\begin{align*}
(\tilde{R}(X,Y)Z)_o &=\tfrac 14 [X,[Y,Z]_{\mathfrak{m}}]_{\mathfrak{m}}-\tfrac 14 [Y,[X,Z]_{\mathfrak{m}}]_{\mathfrak{m}}\\
&-\tfrac 12 [[X,Y]_{\mathfrak{m}},Z]_{\mathfrak{m}}-[[X,Y]_{\mathfrak{h}},Z], \qquad X, Y, Z \in \mathfrak{m}.
\end{align*}
It now follows straightforwardly, by replacing $X,Y,Z$ by $m_1,\dots,m_6$, that using the previously introduced structures the curvature tensor can be written as:
\begin{equation} \label{curvtensorflag}
\begin{split}
\bar R(X,Y)Z&=\tfrac 14 (g(Y,Z)X-g(X,Z)Y)\\
&-\tfrac 14 (g(JY,Z)JX-g(JX,Z)JY +2 g(X,JY)JZ)\\
&+\tfrac 12 (g(J_1 Y,Z)J_1X-g(J_1X,Z)J_1Y +2 g(X,J_1Y)J_1Z)\\
&+\tfrac 12 (g(J_2 Y,Z)J_2X-g(J_2X,Z)J_2Y +2 g(X,J_2Y)J_2Z)\\
&+\tfrac 12 (g(J_3 Y,Z)J_3X-g(J_3X,Z)J_3Y +2 g(X,J_3Y)J_3Z)
\end{split}
\end{equation}

To conclude this section we will investigate how $J_1$, $J_2$ and $J_3$ are related to the Levi-Civita connection of the nearly K\"ahler metric and the skew-symmetric tensor field $G$. First note that as $J=-J_1-J_2-J_3$, we have from the properties of $G$ and $J$ that:
\begin{equation}
G(X,J_1Y)+G(X,J_2Y) + G(X,J_3 Y)=G(J_1 X,Y)+G(J_2 X,Y)+G(J_3 X,Y).
\end{equation}

By a straightforward computation, verifying all components, we get that 
\begin{align} \label{compatibilityJ1G}
J_1 G(X,Y)&=G(J_1 X,Y)+G(X,J_1 Y)+G(X,JY)\\
J_2 G(X,Y)&=G(J_2 X,Y)+G(X,J_2 Y)+G(X,JY)\\
J_3 G(X,Y)&=G(J_3 X,Y)+G(X,J_3 Y)+G(X,JY)
\end{align}
and
\begin{align} \label{nablaJ1}
(\bar \nabla_X J_1)Y&=-\tfrac 12 G(X,Y)-\tfrac 12 JG(J_1X,Y)\\
(\bar \nabla_X J_2)Y&=-\tfrac 12 G(X,Y)-\tfrac 12 JG(J_2X,Y)\\
(\bar \nabla_X J_3)Y&=-\tfrac 12 G(X,Y)-\tfrac 12 JG(J_3X,Y)
\end{align}

\begin{remark} The pseudo-Riemannian counterpart is constructed in a similar way from $SU(2,1)$.  The only differences are that now we start with as a basis of $\mathfrak{su(2,1)}$ the following vector fields:
\begin{alignat*}{2}
&h_1=\begin{pmatrix}
-i&0&0\\0&0&0\\0&0&i
\end{pmatrix}\qquad
&h_2=\begin{pmatrix}
\tfrac{i}{\sqrt{3}}&0&0\\0&\tfrac{-2i}{\sqrt{3}}&0\\0&0&\tfrac{i}{\sqrt{3}}
\end{pmatrix}\\
&m_1=\begin{pmatrix}
0&-1&0\\1&0&0\\0&0&0
\end{pmatrix}\qquad
&m_2=\begin{pmatrix}
0&0&0\\0&0&1\\0&1&0
\end{pmatrix}\\
&m_3=\begin{pmatrix}
0&0&1\\0&0&0\\1&0&0
\end{pmatrix}\qquad
&m_4=\begin{pmatrix}
0&i&0\\i&0&0\\0&0&0
\end{pmatrix}\\
&m_5=\begin{pmatrix}
0&0&0\\0&0&-i\\0&i&0
\end{pmatrix}\qquad
&m_6=\begin{pmatrix}
0&0&-i\\0&0&0\\i&0&0
\end{pmatrix}
\end{alignat*}
Furthermore, the metric must be changed to $g_p (X,Y) = \frac{1}{2} \text{trace} ( I_{-} \overline{X}^T I_{-} Y )$ with $I_{-} = \text{diag} (+1, +1, -1)$ in correspondence with the definition of $\mathfrak{su} (2,1)$. The resulting nearly K\"ahler metric has a signature $4$. In particular it is positive definite on the distribution $V_1$ and negative definite on the distributions $V_2$ and $V_3$. 
\end{remark}

\section{Totally geodesic almost complex surfaces}

Let $M$ be an almost complex surface of the flag manifold. In what follows, we will introduce an $\epsilon$-notation. For the Riemannian flag manifold we will let $\epsilon =1$, whereas for the indefinite flag manifold we will take $\epsilon =-1$. Note that as $g(JX,JX)=g(X,X)$ and $g(X,JX)=0$, in the pseudo-Riemannian case (because we assume that the immersion is a pseudo-Riemannian immersion) the induced metric on $M$ is either positive definite or negative definite. As $J$ preserves the tangent space it follows immediately that
$$h(X,JX)=J h(X,X),$$
where $V$ is a tangent vector to $M$ and $h$ denotes the second fundamental form. Therefore it follows from the Gauss equation that
\begin{equation}\label{secfundform}
K= \langle \tilde{ R}(X,JX)JX,X \rangle -2 \Vert h(X,X ) \Vert^2,
\end{equation}
where $X$ is a unit tangent vector to $M$. Therefore we remark that an almost complex surface is totally geodesic if and only if $K= \langle \tilde{ R}(X,JX)JX,X \rangle$. 

Let $X$ be a unit tangent vector. Note that in the pseudo-Riemannian case when the induced metric on $M$ is negative definite we take $X$ with length $-1$. First if necessary by restricting to an open dense subset, we may assume that locally either:
\begin{enumerate}
\item $X$ and therefore also $TM,$ is contained in a single distribution;
\item $X$ and therefore also $TM,$ is contained in the direct sum of the two distributions;
\item $X$ and therefore also $TM,$ is contained in the direct sum of all three the distributions.
\end{enumerate}
However as we will show that the norm of the orthogonal projections of $X$ onto  the different distributions is constant, this is open dense subset has to be $M$. 

In order to treat all the cases at the same time we write
$$X=a Y +b Z +c W,$$ 
where $Y \in V_1$, $Z \in V_2$ and $W \in V_3$. Moreover we may assume that $Y$, $Z$ and $W$ have unit length. 

Furthermore, we have \begin{equation} \label{eq:ke20} \begin{split} J_1 JX &= - a Y + b Z + c W \\ J_2 JX &= a Y + b Z - c W \\ J_3 JX &= a Y - b Z + c W. \end{split} \end{equation}
In addition, we have \begin{equation} \label{eq:ke21} \begin{split} g(X,X) &= a^2 + \epsilon b^2 + \epsilon c^2 \\ g (X, J_1 J X) &= - a^2 + \epsilon b^2 + \epsilon c^2 \\ g(X, J_2 JX) &= a^2 + \epsilon b^2 - \epsilon c^2 \\ g(X, J_3 J X) &= a^2 - \epsilon b^2 + \epsilon c^2. \end{split} \end{equation}
The Riemann tensor becomes \begin{equation} \label{eq:ke22} \begin{split} R(X, JX) JX &= - \frac{1}{2} (a^2 + \epsilon b^2 + \epsilon c^2) X + \frac{3}{2} \left( 3 a^3 - \epsilon a b^2 - \epsilon a c^2 \right) Y \\ & + \frac{3}{2} \left(  3 \epsilon b^3 - a^2 b - \epsilon b c^2 \right) Z + \frac{3}{2} \left( 3 \epsilon c^3 - \epsilon b^2 c - a^2 c \right) W. \end{split} \end{equation}
As $M$ is totally geodesic, we must have that $R(X,JX)JX$ is tangential and therefore because of the properties of the curvature tensor, we need that the matrix \begin{equation} \label{eq:ke23} \begin{pmatrix} \frac{3}{2} \left( 3 a^3 - \epsilon a b^2 - \epsilon a c^2 \right) & \frac{3}{2} \left( 3 \epsilon  b^3 - a^2 b - \epsilon b c^2 \right) & \frac{3}{2} \left( 3 \epsilon c^3 - \epsilon b^2 c - a^2 c \right) \\ a & b & c \end{pmatrix} \end{equation} to have rank $1$. This happens exactly when the determinant of all of the $(2 \times 2$)-minors vanish. This leads to the equations
\begin{align}
 &  a(a^2 - \epsilon b^2) b = 0 \label{eq:i}
   \\
  & a (a^2  - \epsilon c^2) c  = 0 \label{eq:ii} \\
  & b (c^2 - b^2) c= 0. \label{eq:iii} 
\end{align}
Note how the last equation is independent of $\epsilon$. This corresponds to the case when $T_p M$ lies in the last two distributions, while the first two correspond to when $T_p M \subset V_1 (p) \oplus V_2 (p)$ or $T_p M \subset V_1 (p) \oplus V_3(p)$. Indeed: the metric on the last two distributions is either positive-definite in the Riemannian case $( \epsilon = +1)$ or negative definite in the pseudo-Riemannian case $(\epsilon = -1)$, so essentially this is just a change in an overall minus sign. However, when $T_p M$ lies in the first two distributions for example, there is a relative minus sign difference in going from the Riemannian to the pseudo-Riemannian case. Also, when $T_p M$ lies in two distributions of different signature, we do not expect totally geodesic almost complex surfaces to exist. More specifically, we have: 

\begin{enumerate}
\item If $(a,b,c) = (1,0,0)$ or any other permutation, then $K  = 4$ as the reader may check.
\item Assume that $c = 0$ but $a,b \neq 0$. Then we need to to satisfy equation \eqref{eq:i} above, which can only be done when $\epsilon = 1$, i.e. in the Riemannian case. In the Riemannian case, we can then pick $a = b = 1 / \sqrt{2}$ and find $K = 1$. Similarly for $b = 0$ but $a, c \neq 0$. Finally, if $a = 0$ but $b, c \neq 0$, then we need to satisfy equation \eqref{eq:iii} above, which can be done in both the Riemannian and semi-Riemannian case. Then we again find $K = 1$. 
\item Finally assume that $a, b, c \neq 0$, so that $T_p M$ lies in all three distributions. It is clear that only all three equations can be satisfied simultaneously when $\epsilon = +1$, i.e. in the Riemannian case, with $(a,b,c) = ( \frac{1}{ \sqrt{3}}, \frac{1}{ \sqrt{3}}, \frac{1}{ \sqrt{3}})$. In that case, we find $K = 0$ as the reader can check. In the semi-Riemannian case it is impossible to satisfy all three simultaneously. 
\end{enumerate} 

Note that as the flag manifold is homogeneous we may of course assume that the totally geodesic surface contains the point $\pi(I)$. Then summarizing the previous arguments in the Riemannian case we have either:
\begin{enumerate}
\item $TM \subset V_1$, the tangent space $T$ at $\pi(I)$ is spanned by $m_1$ and $m_4$, and $M$ is locally isometric with a sphere of curvature $4$ (the other cases where the tangent space is contained in one distribution are equivalent);
\item $TM \subset V_1 \oplus V_2$, by applying an isometry (which preserves the identity), corresponding to a rotation in the first two distributions, we may assume that the tangent space $T$ at the identity is spanned by $\frac{1}{ \sqrt{2}}(m_1+m_2)$ and $\frac{1}{ \sqrt{2}}(m_4+m_5)$. Moreover  $M$ is locally isometric with a sphere of curvature $1$ (the other cases where the tangent space is contained in different sums of two distributions are equivalent);
\item $TM \subset V_1 \oplus V_2 \oplus V_3$. Recall that
$X=\tfrac{1}{\sqrt{3}} (Y+Z+W)$, where $Y,Z,W$ are unit vectors belonging to the appropriate distributions. As $G(m_1,m_2)=  m_6$, it follows that we can write
$G(Y,Z)=\cos \theta W + \sin \theta  JW,$
for some function $\theta$. We now define $\tilde X= \cos \varphi X + \sin \varphi JX$.
It then follows that
\begin{align*}
&\tilde Y= \cos \varphi Y + \sin \varphi JY\\
&\tilde Z= \cos \varphi Z + \sin \varphi JZ\\
&\tilde W= \cos \varphi W + \sin \varphi JW\\
&G(\tilde Y,\tilde Z)=\cos (2\varphi-\theta) W +\sin(2\varphi-\theta) JW\\
&G(\tilde Y,\tilde Z)=\cos (\varphi-\theta) \tilde W +\sin(\varphi-\theta) J\tilde W
\end{align*}
This shows that we can choose $X$ such that $G(Y,Z)=JW$.
So by applying an isometry (which preserves the identity), corresponding to a rotation in the first two distributions, we may assume that the tangent space $T$ at the identity is spannned by $\frac{1}{ \sqrt{3}}(m_1+m_2+m_3)$ and $\frac{1}{ \sqrt{3}}(m_4+m_5-m_6)$, and $M$ is locally equivalent to a flat space. 
\end{enumerate}

In a similar way in the semi-Riemannian case we are left with the following possibilities:
\begin{enumerate}
\item $TM \subset V_1$, the tangent space $T$ at $\pi(I)$ is spanned by $m_1$ and $m_4$, and $M$ is locally isometric with a sphere of curvature $4$,
\item $TM \subset V_2$, the tangent space $T$ at $\pi(I)$ is spanned by $m_2$ and $m_5$, and $M$ is locally anti-isometric with a hyperbolic plane of curvature $-4$. The case that $TM \subset V_3$
leads to an equivalent result. 
\item $TM \subset V_2 \oplus V_3$, by applying an isometry (which preserves the identity), corresponding to a rotation in the first two distributions, we may assume that the tangent space $T$ at the identity is spannned by $\frac{1}{ \sqrt{2}}(m_2+m_3)$ and $\frac{1}{ \sqrt{2}}(m_5-m_6)$. Moreover  $M$ is locally anti-isometric with a hyperbolic plane of curvature $-1$.
\end{enumerate}

Note that as a (pseudo)-Riemannian submersion maps horizontal geodesics to geodesics of the flag manifold, in order to complete the classification, it is now sufficient to check that the surface obtained by taking $\pi (\text{exp}_{I}(t T))$ is indeed a totally geodesic almost complex surface. By construction and the Hopf-Rinow theorem, as the induced metric on the surface is definite and the geodesics through the point $\pi(I)$ are defined for all time, these surfaces are complete. Moreover, as $SU(3)$ in the Riemannian case (and $SU(2,1)$ in the pseudo Riemannian case) acts by isometries on the flag manifold preserving the nearly K\"ahler structure, by construction it is also clear that the image is homogeneous and almost complex. 

We will now give explicit parametrizations of these surfaces and verify that the immersions are indeed totally geodesic and almost complex. 
 
\begin{example}
We look at the surface $F_1(t,u)=\exp(t(\cos(u) m_1+\sin(u) m_4))$ in $SU(3)$. It follows that
$$F_1(t,u)=\left(
\begin{array}{ccc}
 \cos (t) & -e^{-i u} \sin (t) & 0 \\
 e^{i u} \sin (t) & \cos (t) & 0 \\
 0 & 0 & 1 \\
\end{array}
\right)$$
It follows that
\begin{align*}
&dF_1 \left( \frac{\partial}{\partial t} \right)= F_1(t,u) \cdot \left(
\begin{array}{ccc}
 0 & -e^{-i u} & 0 \\
 e^{i u} & 0 & 0 \\
 0 & 0 & 0 \\
\end{array}
\right)\\
&dF_1 \left( \frac{\partial}{\partial u} \right)=F_1(t,u) \cdot \left(
\begin{array}{ccc}
 i \sin ^2(t) & \sin (t) \cos (t) ie^{-i u} & 0 \\
 i e^{i u} \sin (t) \cos (t) & -i \sin ^2(t) & 0 \\
 0 & 0 & 0 \\
\end{array}
\right).
\end{align*}
Note that $dF_1(\frac{\partial}{\partial t})$ is horizontal and that the horizontal part $(dF_1(\frac{\partial}{\partial u}))^h$ is given by
$$ \left( dF_1 \left( \frac{\partial}{\partial u} \right) \right)^h=F_1(t,u) \cdot \left(
\begin{array}{ccc}
 0& \sin (t) \cos (t) ie^{-i u} & 0 \\
 i e^{i u} \sin (t) \cos (t) & 0 & 0 \\
 0 & 0 & 0 \\
\end{array}\right). $$
From this we immediately see that the surface $\pi(F_1(t,u))$ is indeed an almost complex surface from which the tangent space always lies in the first distribution. We also have that the induced metric on the surface is given by
\begin{align*}
& < \frac{\partial}{\partial t},\frac{\partial}{\partial t} > =1\\
& < \frac{\partial}{\partial t},\frac{\partial}{\partial u} > =0\\
& < \frac{\partial}{\partial u},\frac{\partial}{\partial u} > =\left( \tfrac{ \sin(2t)}{2} \right)^2,
\end{align*}
which has constant curvature $4$. Applying now \eqref{secfundform} shows that the surface is indeed totally geodesic.
\end{example}

\begin{example}
We look at the surface $F_2(t,u)=\exp(t(\cos(u) \tfrac{1}{\sqrt{2}}(m_1+m_2)+ \sin(u)  \tfrac{1}{\sqrt{2}}(m_4+m_5)))$ in $SU(3)$. It follows that
$$F_2(t,u)=\left(
\begin{array}{ccc}
 \cos ^2\left(\frac{t}{2}\right) & -\frac{e^{-i u} \sin (t)}{\sqrt{2}} & e^{-2 i u} \sin
   ^2\left(\frac{t}{2}\right) \\
 \frac{e^{i u} \sin (t)}{\sqrt{2}} & \cos (t) & -\frac{e^{-i u} \sin (t)}{\sqrt{2}} \\
 e^{2 i u} \sin ^2\left(\frac{t}{2}\right) & \frac{e^{i u} \sin (t)}{\sqrt{2}} & \cos
   ^2\left(\frac{t}{2}\right) \\
\end{array}
\right)$$
One calculates that
\begin{align*}
&dF_2 \left( \frac{\partial}{\partial t} \right)= F_2(t,u) \cdot \left(
\begin{array}{ccc}
 0 & -\frac{e^{-i u}}{\sqrt{2}} & 0 \\
 \frac{e^{i u}}{\sqrt{2}} & 0 & -\frac{e^{-i u}}{\sqrt{2}} \\
 0 & \frac{e^{i u}}{\sqrt{2}} & 0 \\
\end{array}
\right)\\
&dF_2 \left( \frac{\partial}{\partial u} \right) =F_2(t,u) \cdot \left(
\begin{array}{ccc}
 -i (\cos (t)-1) & \frac{i e^{-i u} \sin (t)}{\sqrt{2}} & 0 \\
 \frac{i e^{i u} \sin (t)}{\sqrt{2}} & 0 & \frac{i e^{-i u} \sin (t)}{\sqrt{2}}
   \\
 0 & \frac{i e^{i u} \sin (t)}{\sqrt{2}} & i (\cos (t)-1) \\
\end{array}
\right)
\end{align*}
Note that $dF_2(\frac{\partial}{\partial t})$ is horizontal and that the horizontal part $(dF_2(\frac{\partial}{\partial u}))^h$ is given by
$$ \left( dF_2 \left( \frac{\partial}{\partial u} \right) \right)^h=F_2(t,u) \cdot \left(
\begin{array}{ccc}
0 & \frac{i e^{-i u} \sin (t)}{\sqrt{2}} & 0 \\
 \frac{i e^{i u} \sin (t)}{\sqrt{2}} & 0 & \frac{i e^{-i u} \sin (t)}{\sqrt{2}}
   \\
 0 & \frac{i e^{i u} \sin (t)}{\sqrt{2}} &0 \\
\end{array}
\right)$$
From this we immediately see that the surface $\pi(F_2(t,u))$ is indeed an almost complex surface from which the tangent space always lies in the direct sum of the first two distribution in the desired way. We also have that the induced metric on the surface is given by
\begin{align*}
&<\frac{\partial}{\partial t},\frac{\partial}{\partial t}>=1\\
&<\frac{\partial}{\partial t},\frac{\partial}{\partial u}>=0\\
&<\frac{\partial}{\partial u},\frac{\partial}{\partial u}>=\sin^2 t ,
\end{align*}
which has constant curvature $1$. Applying now \eqref{secfundform} shows that the surface is indeed totally geodesic.
\end{example}

\begin{example}
We look at the surface $F_3(t,u)= \exp(t \tfrac{1}{\sqrt{3}}(m_1+m_2+m_3)+ u  \tfrac{1}{\sqrt{3}}(m_4+m_5-m_6))$ in $SU(3)$. It follows that $F_3(t,u)=(c_{ij})$, where
\begin{align*}
c_{11}&=
 \frac{1}{3} e^{-\frac{2 i u}{\sqrt{3}}} \left(1+2 e^{i \sqrt{3} u} \cos (t)\right)\\
c_{12}&=
   \frac{1}{3} e^{-\frac{2 i u}{\sqrt{3}}} \left(-1+e^{i \sqrt{3} u} \left(\cos (t)-\sqrt{3}
   \sin (t)\right)\right)\\
c_{13}    &= \frac{1}{3} e^{-i t-\frac{2 i u}{\sqrt{3}}} \left(-e^{i
   \left(t+\sqrt{3} u\right)} \left(\sqrt{3} \sin (t)+\cos (t)\right)+i \sin (t)+\cos
   (t)\right) \\
c_{21}&=  
 \frac{1}{3} e^{-\frac{2 i u}{\sqrt{3}}} \left(-1+e^{i \sqrt{3} u} \left(\sqrt{3} \sin
   (t)+\cos (t)\right)\right)\\
c_{22}   &= \frac{1}{3} e^{-\frac{2 i u}{\sqrt{3}}} \left(1+2 e^{i
   \sqrt{3} u} \cos (t)\right) \\
 c_{23}  &= \frac{1}{3} e^{-\frac{2 i u}{\sqrt{3}}} \left(-1+e^{i
   \sqrt{3} u} \left(\cos (t)-\sqrt{3} \sin (t)\right)\right) \\
c_{31}&= \frac{1}{3} e^{-i t-\frac{2 i u}{\sqrt{3}}} \left(e^{i \left(t+\sqrt{3} u\right)}
   \left(\sqrt{3} \sin (t)-\cos (t)\right)+i \sin (t)+\cos (t)\right)\\
  c_{32}&=   \frac{1}{3}
   e^{-\frac{2 i u}{\sqrt{3}}} \left(-1+e^{i \sqrt{3} u} \left(\sqrt{3} \sin (t)+\cos
   (t)\right)\right)\\
c_{33}    &= \frac{1}{3} e^{-\frac{2 i u}{\sqrt{3}}} \left(1+2 e^{i \sqrt{3} u}
   \cos (t)\right)
\end{align*}
It follows that
\begin{align*}
&dF_3 \left( \frac{\partial}{\partial t} \right) = F_3(t,u) \cdot \left(
\begin{array}{ccc}
 0 & -\frac{1}{\sqrt{3}} & -\frac{1}{\sqrt{3}} \\
 \frac{1}{\sqrt{3}} & 0 & -\frac{1}{\sqrt{3}} \\
 \frac{1}{\sqrt{3}} & \frac{1}{\sqrt{3}} & 0 \\
\end{array}
\right)\\
&dF_3 \left( \frac{\partial}{\partial u} \right) =F_3(t,u) \cdot \left(
\begin{array}{ccc}
 0 & \frac{i}{\sqrt{3}} & -\frac{i}{\sqrt{3}} \\
 \frac{i}{\sqrt{3}} & 0 & \frac{i}{\sqrt{3}} \\
 -\frac{i}{\sqrt{3}} & \frac{i}{\sqrt{3}} & 0 \\
\end{array}
\right)
\end{align*}
Note that $dF_3(\frac{\partial}{\partial t})$  and   $dF_3(\frac{\partial}{\partial u})$ are horizontal in this case. 
It is also immediately clear  that the surface $\pi(F_3(t,u))$ is indeed an almost complex surface from which the tangent space always lies in the direct sum of the three distributions in the desired way. We also have that the induced metric on the surface is given by
\begin{align*}
&<\frac{\partial}{\partial t},\frac{\partial}{\partial t}>=1\\
&<\frac{\partial}{\partial t},\frac{\partial}{\partial u}>=0\\
&<\frac{\partial}{\partial u},\frac{\partial}{\partial u}>=1,
\end{align*}
which is a flat metric. It is also clear from the expression of $F_3$ that the image is a torus. Applying now \eqref{secfundform} shows that the surface is indeed totally geodesic. 
\end{example}

Using similar computations, in the pseudo Riemannian case we obtain the following examples.
\begin{example}
We look at the surface $F_4(t,u)= \exp(t(\cos(u) m_1+\sin(u) m_4)) $ in $SU(2,1)$. It follows that
$$F_4(t,u)=\left(
\begin{array}{ccc}
 \cos (t) & -e^{-i u} \sin (t) & 0 \\
 e^{i u} \sin (t) & \cos (t) & 0 \\
 0 & 0 & 1 \\
\end{array}
\right)$$
It follows that
\begin{align*}
&dF_4 \left( \frac{\partial}{\partial t} \right) = F_4(t,u) \cdot \left(
\begin{array}{ccc}
 0 & -e^{-i u} & 0 \\
 e^{i u} & 0 & 0 \\
 0 & 0 & 0 \\
\end{array}
\right)\\
&dF_4 \left( \frac{\partial}{\partial u} \right) =F_4(t,u) \cdot \left(
\begin{array}{ccc}
 i \sin ^2(t) & \sin (t) \cos (t) ie^{-i u} & 0 \\
 i e^{i u} \sin (t) \cos (t) & -i \sin ^2(t) & 0 \\
 0 & 0 & 0 \\
\end{array}
\right)
\end{align*}
Note that $dF_4(\frac{\partial}{\partial t})$ is horizontal and that the horizontal part $(dF_4(\frac{\partial}{\partial u}))^h$ is given by
$$ \left( dF_4 \left( \frac{\partial}{\partial u} \right) \right)^h=F_4(t,u) \cdot \left(
\begin{array}{ccc}
 0& \sin (t) \cos (t) ie^{-i u} & 0 \\
 i e^{i u} \sin (t) \cos (t) & 0 & 0 \\
 0 & 0 & 0 \\
\end{array}\right)$$
From this we immediately see that the surface $\pi(F_4(t,u))$ is indeed an almost complex surface from which the tangent space always lies in the first distribution. We also have that the induced metric on the surface is given by
\begin{align*}
&<\frac{\partial}{\partial t},\frac{\partial}{\partial t}>=1\\
&<\frac{\partial}{\partial t},\frac{\partial}{\partial u}>=0\\
&<\frac{\partial}{\partial u},\frac{\partial}{\partial u}>= \left( \tfrac{ \sin(2t)}{2} \right)^2,
\end{align*}
which has constant curvature $4$. Applying now \eqref{secfundform} shows that the surface is indeed totally geodesic.
\end{example}

\begin{example}
We look at the surface $F_5(t,u)= \exp(t(\cos(u) m_2+\sin(u) m_5)$ in $SU(2,1)$. It follows that
$$F_5(t,u)=\left(
\begin{array}{ccc}
 1 & 0 & 0 \\
 0 & \cosh (t) & e^{-i u} \sinh (t) \\
 0 & e^{i u} \sinh (t) & \cosh (t) \\
\end{array}
\right)$$
One calculates that
\begin{align*}
&dF_5 \left( \frac{\partial}{\partial t} \right)= F_5(t,u) \cdot \left(
\begin{array}{ccc}
 0 & 0 & 0 \\
 0 & 0 & e^{-i u} \\
 0 & e^{i u} & 0 \\
\end{array}
\right)\\
&dF_5 \left( \frac{\partial}{\partial u} \right)=F_5(t,u) \cdot \left(
\begin{array}{ccc}
 0 & 0 & 0 \\
 0 & -i \sinh ^2(t) & -i e^{-i u} \sinh (t) \cosh (t) \\
 0 & i e^{i u} \sinh (t) \cosh (t) & i \sinh ^2(t) \\
\end{array}
\right)
\end{align*}
Note that $dF_5(\frac{\partial}{\partial t})$ is horizontal and that the horizontal part $(dF_5(\frac{\partial}{\partial u}))^h$ is given by
$$ \left( dF_5 \left( \frac{\partial}{\partial u} \right) \right)^h=F_5(t,u) \cdot \left(
\begin{array}{ccc}
 0 & 0 & 0 \\
 0 & 0 & -i e^{-i u} \sinh (t) \cosh (t) \\
 0 & i e^{i u} \sinh (t) \cosh (t) & 0 \\
\end{array}
\right)$$
From this we immediately see that the surface $\pi(F_5(t,u))$ is indeed an almost complex surface from which the tangent space always lies in the second distribution. We also have that the induced metric on the surface is given by
\begin{align*}
&<\frac{\partial}{\partial t},\frac{\partial}{\partial t}>=-1\\
&<\frac{\partial}{\partial t},\frac{\partial}{\partial u}>=0\\
&<\frac{\partial}{\partial u},\frac{\partial}{\partial u}>=- \left( \tfrac{ \sinh(2t)}{2} \right)^2,
\end{align*}
which has constant curvature $4$. As the metric is negative definite it follows that the surface is anti-isometric with a hyperbolic plane. Applying now \eqref{secfundform} shows that the surface is indeed totally geodesic.
\end{example}

\begin{example}
We look at the surface $F_6(t,u)= \exp(t(\cos(u) \tfrac{1}{\sqrt{2}}(m_2+m_3)+ \sin(u)  \tfrac{1}{\sqrt{2}}(m_5-m_6))) $ in $SU(2,1)$. It follows that
$$F_6(t,u)=\left(
\begin{array}{ccc}
 \cosh ^2\left(\frac{t}{2}\right) & e^{2 i u} \sinh ^2\left(\frac{t}{2}\right) & \frac{e^{i
   u} \sinh (t)}{\sqrt{2}} \\
 e^{-2 i u} \sinh ^2\left(\frac{t}{2}\right) & \cosh ^2\left(\frac{t}{2}\right) &
   \frac{e^{-i u} \sinh (t)}{\sqrt{2}} \\
 \frac{e^{-i u} \sinh (t)}{\sqrt{2}} & \frac{e^{i u} \sinh (t)}{\sqrt{2}} & \cosh (t) \\
\end{array}
\right)$$
It follows that
\begin{align*}
&dF_6 \left( \frac{\partial}{\partial t} \right) = F_6(t,u) \cdot \left(
\begin{array}{ccc}
 0 & 0 & \frac{e^{i u}}{\sqrt{2}} \\
 0 & 0 & \frac{e^{-i u}}{\sqrt{2}} \\
 \frac{e^{-i u}}{\sqrt{2}} & \frac{e^{i u}}{\sqrt{2}} & 0 \\
\end{array}
\right)\\
&dF_6 \left( \frac{\partial}{\partial u} \right)=F_6(t,u) \cdot \left(
\begin{array}{ccc}
 i (\cosh (t)-1) & 0 & \frac{i e^{i u} \sinh (t)}{\sqrt{2}} \\
 0 & -i (\cosh (t)-1) & -\frac{i e^{-i u} \sinh (t)}{\sqrt{2}} \\
 -\frac{i e^{-i u} \sinh (t)}{\sqrt{2}} & \frac{i e^{i u} \sinh (t)}{\sqrt{2}} & 0 \\
\end{array}
\right)
\end{align*}
Note that $dF_6(\frac{\partial}{\partial t})$ is horizontal and that the horizontal part $(dF_6(\frac{\partial}{\partial u}))^h$ is given by
$$ \left( dF_6 \left( \frac{\partial}{\partial u} \right) \right)^h=F_6(t,u) \cdot \left(
\begin{array}{ccc}
0 & 0 & \frac{i e^{i u} \sinh (t)}{\sqrt{2}} \\
 0 & 0 & -\frac{i e^{-i u} \sinh (t)}{\sqrt{2}} \\
 -\frac{i e^{-i u} \sinh (t)}{\sqrt{2}} & \frac{i e^{i u} \sinh (t)}{\sqrt{2}} & 0 \\
\end{array}
\right)$$
From this we immediately see that the surface $\pi(F_6(t,u))$ is indeed an almost complex surface from which the tangent space always lies in the direct sum of the last two distribution in the desired way. We also have that the induced metric on the surface is given by
\begin{align*}
&<\frac{\partial}{\partial t},\frac{\partial}{\partial t}>=-1\\
&<\frac{\partial}{\partial t},\frac{\partial}{\partial u}>=0\\
&<\frac{\partial}{\partial u},\frac{\partial}{\partial u}>=-\sinh^2 t ,
\end{align*}
which has constant curvature $1$. As the induced metric is negative definite this implies that the surface is anti isometric with a hyperbolic plane. Applying now \eqref{secfundform} shows that the surface is indeed totally geodesic.
\end{example}


\begin{thebibliography}{00}

\bibitem{apostolovgrantcharovivanov} V. Apostolov, G. Grantcharov, S.  Ivanov.  Orthogonal complex structures on certain Riemannian 6-manifolds. Differential Geom. Appl.  11  (1999),  no. 3, 279--296.

\bibitem{arvanitoyeorgos} A. Arvanitoyeorgos. An Introduction to Lie groups and the Geometry of Homogeneous Spaces. Student Mathematical Library, Vol. 22, 2003.

\bibitem{besse} A.L. Besse.  Einstein manifolds. Reprint of the 1987 edition. Classics in Mathematics. Springer-Verlag, Berlin,  2008. xii+516 pp. 

\bibitem{butruille} J. B. Butruille. Classification des vari\'{e}t\'{e} approximativement k\"ahleriennes homog\'{e}nes [[Classification of the Nearly-k\"ahler Homogeneties]].  Ann. Global Anal. Geom.  27  (2005),  no. 3, 201--225.

\bibitem{loubeau} G. Deschamps and Eric Loubeau, Hypersurfaces of the nearly K\"ahler twistor spaces $\mathbb CP^3$ and $\mathbb F_{12}$. T\^ohoku Mathematical Journal. To Appear

\bibitem{foshas} L. Foscolo and M.  Haskins. New G2-holonomy cones and exotic nearly K\"ahler structures on $S^6$ and $S^3 \times S^3$.  Ann. of Math. (2)  185  (2017),  no. 1, 59--130.

\bibitem{gray1} A. Gray and Luis M.Hervella, The sixteen classes of almost Hermitian manifolds and their linear invariants.
 Ann. Mat. Pura Appl. (4)  123  (1980), 35--58.
		
\bibitem{gray2} A.  Gray. The structure of nearly K\"ahler manifolds. Math. Ann.  223  (1976),  no. 3, 233--248.
		
\bibitem{gray3} A.  Gray. Riemannian manifolds with geodesic symmetries of order 3. J. Differential Geometry  7  (1972), 343--369.
		
\bibitem{gray4} A. Gray, Alfred. Nearly K\"ahler manifolds. J. Differential Geometry  4  (1970), 283--309.

\bibitem{kobayashinomizu} S. Kobayashi, Shoshichi and K.  Nomizu, Foundations of differential geometry. Vol I and  Vol. II.
Wiley Classics Library. A Wiley-Interscience Publication.
John Wiley \& Sons, Inc., New York,  1996. 


\bibitem{reinierstorm} R. Storm. Lagrangian submanifolds of the nearly K\"ahler full flag manifold $F_{1,2}(\mathbb C^3)$.
 J. Geom. Phys.  158  (2020), 103844, 16 pp.
 
 
\end{thebibliography}
\end{document}